\documentclass[12pt]{article}

\usepackage{amsmath}
\usepackage{amsfonts}
\usepackage{amssymb}
\usepackage{graphicx, color}
\usepackage{hyperref}
\usepackage{epsfig}
\usepackage{epstopdf} 

\definecolor{darkgreen}{rgb}{0,0.55,0}

\newtheorem{proposition}{Proposition}[section]
\newtheorem{theorem}{Theorem}[section]
\newtheorem{lemma}[theorem]{Lemma}

\newtheorem{remark}[theorem]{Remark}

\DeclareSymbolFont{AMSb}{U}{msb}{m}{n}
\DeclareMathSymbol{\N}{\mathbin}{AMSb}{"4E}
\DeclareMathSymbol{\Z}{\mathbin}{AMSb}{"5A}
\DeclareMathSymbol{\R}{\mathbin}{AMSb}{"52}
\DeclareMathSymbol{\Q}{\mathbin}{AMSb}{"51}
\DeclareMathSymbol{\I}{\mathbin}{AMSb}{"49}

\numberwithin{equation}{section}

\begin{document}

\title{\bf A Generalization of the Sphere Covering Inequality}
\author{{Changfeng Gui\footnote{Department of Mathematics, University of Macau, Macau SAR, P. R. China. E-mail address: changfenggui@um.edu.mo}
\qquad Amir Moradifam\footnote{Department of Mathematics, University of California, Riverside, California, USA. E-mail: moradifam@math.ucr.edu. }}}
\date{\today}

\smallbreak \maketitle

\begin{abstract}
The Sphere Covering Inequality was introduced in \cite{GM} (\emph{Invent. Math.}, 2018) as a sharp geometric inequality that provides a lower bound for the total area of two distinct surfaces of Gaussian curvature 1. These surfaces are assumed to be conformal to the Euclidean unit disk and share the same conformal factor along the boundary. In this paper, we establish a quantitative generalization that relaxes the boundary matching condition by allowing the conformal factors to differ by a constant \( c \ge 0 \) on the boundary. This refinement reveals a new stability-type structure underlying the inequality. Our results show that the Sphere Covering Inequality is stable with respect to perturbations in the boundary data and provide a precise quantitative description of how the total-area bound varies under such perturbations. The generalized inequality provides new analytic and geometric tools for the study of elliptic equations with exponential nonlinearities, conformal geometry, and related problems in mathematical physics.

\end{abstract}

\section{Introduction}

Second-order nonlinear elliptic equations with exponential nonlinearities form one of the most fundamental classes of equations in two-dimensional analysis. Equations of the form  
\begin{equation}\label{firstPDE}
\Delta u + e^{u} = f(x), \quad x \in \Omega \subset \mathbb{R}^2,
\end{equation}
arise in a wide variety of problems in mathematics, mathematical physics, and biology.  Equations of this form naturally arise in problems where curvature, conformal geometry, and analytic inequalities interact. In differential geometry, they describe surfaces with prescribed Gaussian curvature.  Analytically, they are closely related to the Moser–Trudinger type inequalities that capture the critical exponential growth characterizing the limiting case of Sobolev embeddings in dimension two. The interplay between geometric structures and sharp functional inequalities in this setting has led to deep insights into the nature of conformal invariance and the structure of extremal functions, see \cite{Aubin-MR0448404, Aubin2-MR534672, Beckner-MR1230930,  ChangGui-CPA2023, CL1-MR1121147, CL2-MR1338474, CLiu-MR1209959, CY1-MR908146, CY2-MR925123, DEJ-MR3377875, DET-MR2509375, GL-MR2670931, GW-MR1760786, Moser-MR0301504, O-MR677001, OPS-MR960228, T-MR0216286}.

In mathematical physics, such equations model a wide range of nonlinear phenomena. They govern the Chern–Simons self-dual vortex equations, Liouville and Toda systems, and mean-field equations arising in the statistical mechanics of two-dimensional turbulence and vortex dynamics. They also appear in the theory of self-gravitating cosmic strings in Einstein’s general relativity, where the exponential nonlinearity encodes the coupling between geometry and matter fields. Beyond these physical contexts, similar structures emerge in the theory of elliptic and hyperelliptic functions and in free-boundary problems such as models of cell motility and biological pattern formation; see \cite{BartolucciChenLinTarantello-MR2071339, BartolucciLin-MR2542689, BartolucciTarantello-MR1917679, BD-MR3355004, BerlyandFuhrmannRybalko-MR3838409, BLT-MR2838340, CabreLuciaSanchon-MR2153523, CaffarelliYang-MR1328259, CCL-MR2055839, CGS-MR2944090, ChanFuLin-MR1946446, CK-MR1262195, CLMP1-MR1145596, CLMP2-MR1362165, DingJostLiWang-MR1634402, K-MR1193342, Lin1-MR1770683, LinLucia-MR2263535, LinWang-MR2680439, LinWang-MR3665624, LinWeiYe-MR2969270, PT-MR2876669, Yang-MR1277473}.

Because of their critical nonlinearity and deep geometric origins, these equations provide a unifying framework that connects diverse areas of analysis, geometry, and mathematical physics. Their study has led to the discovery of sharp functional inequalities, refined classification results, and symmetry and uniqueness theorems for Liouville-type and mean-field equations. The interplay between the analytic structure of \eqref{firstPDE} and its geometric and physical interpretations continues to inspire new methods and results in modern nonlinear analysis.

A fundamental breakthrough in this area was the introduction of the \emph{Sphere Covering Inequality (SCI)} in~\cite{GM}, which asserts that the combined area of two distinct surfaces with Gaussian curvature equal to one, conformal to the Euclidean unit disk, and sharing the same conformal factor on the boundary, must be at least~$4\pi$. Geometrically, this means that the two surfaces together cover the entire unit sphere, up to a conformal rearrangement. Indeed, the following theorem was proved in \cite{GM}. \\

\begin{theorem}[The Sphere Covering Inequality \cite{GM}]\label{MainTheorem0}
Let $\Omega $ be a simply-connected  subset of $\R^2$ and assume $v_i \in C^2(\overline{\Omega})$, $i=1,2$ satisfy
\begin{equation}
\Delta v_i +e^{2v_i}=f_{i}(y),
\end{equation}
where $f_2 \geq f_1 \ge 0$  in $\Omega$.   If $v_2 \ge v_1,  v_2 \not \equiv v_1$ in $\omega$ and $v_2=v_1$ on $\partial \omega$  for some  piecewise
Liptschitz subdomain $\omega \subset \Omega$,  then
\begin{equation}
\int_{\omega} (e^{2v_1}+e^{2v_2})dy \geq 4\pi. 
\end{equation}
Moreover,  the equality only holds when $f_2 \equiv f_1 \equiv 0$  in $\omega$, and $(\omega, e^{2v_i}dy ) $, $i=1,2$ are isometric to  two complementary spherical caps on the standard unit sphere.  
\end{theorem}

The Sphere Covering Inequality was recently introduced in~\cite{GM} and has been applied to address several problems concerning symmetry and uniqueness of solutions to elliptic equations with exponential nonlinearities in dimension \( n = 2 \). In particular, it was employed to resolve a long-standing conjecture of Chang and Yang~\cite{CY2-MR925123} regarding the best constant in Moser--Trudinger type inequalities~\cite{GM}, as well as a conjecture of Bartolucci, Lin, and Tarantello~\cite{BLT-MR2838340} concerning mean field equations with singularities on~\(\mathbb{S}^2\) arising from self-dual gauge field theory. Moreover, it has led to a variety of further symmetry and uniqueness results for mean field equations, Onsager vortex models, the Sinh--Gordon and cosmic string equations, Toda systems, and rigidity phenomena related to the Hawking mass in general relativity (see \cite{BartolucciGuiJevnikarMoradifam-MR3961327, GuiJevnikarMoradifam-MR3769320, GuiMoradifam-MR3762223, GuiMoradifam-MR3927157, GuiMoradifam-MR3741415, LeeLinTarantelloYang-MR3699389, ShiSunTianWei-MRXXXXXXX, WangWangYang-arXiv1804.03319}). Beyond these specific applications, the Sphere Covering Inequality represents a central analytic principle in the treatment of nonlinear PDEs with critical exponential growth in two dimensions.

In this paper, we establish a generalization of the Sphere Covering Inequality that relaxes the boundary condition requiring two conformal metrics (or solutions) to coincide exactly on the boundary. Instead, we allow the solutions to differ by a constant \(c \ge 0\) on \(\partial \Omega\). This generalization broadens the applicability of the SCI to problems with nonidentical boundary data while preserving its geometric sharpness. Such a refinement provides new analytic tools for the study of elliptic PDEs with relaxed boundary conditions, symmetry-breaking phenomena, and geometric stability problems arising in conformal geometry and related fields.

Specifically, instead of assuming \(v_1 = v_2\) on \(\partial \Omega\), we permit a constant boundary gap:
\[
v_2 - v_1 = c \ge 0 \quad \text{on } \partial \Omega.
\]
Our main result shows that under these weaker boundary conditions one still obtains a quantitative lower bound, but with the constant \(4\pi\) replaced by a parameter depending on the boundary gap \(c\) relative to the maximal difference of the solutions inside \(\Omega\). Precisely, if \(\omega \subset \Omega\) is a Lipschitz subdomain on which \(v_2 \ge v_1\) and \(v_2 \not\equiv v_1\), then
\[
\int_{\omega}\big(e^{2v_1}+e^{2v_2}\big)\,dy \;\ge\; 4(1 - \frac{c}{M})\pi,
\]
where
\[
 M = \max_{\omega}(v_2 - v_1).
\]
When \(c = 0\), the result reduces to the classical Sphere Covering Inequality, and if $c=\frac{1}{2}\max\limits_{\omega} (v_2 - v_1)>0,$ then 
\[\int_{\omega} (e^{2v_1}+e^{2v_2})dy \geq 2 \pi.\]
Indeed we prove the following theorem which is the main result of this paper.  \\
\begin{theorem}[Generalized Sphere Covering Inequality]\label{GeneralSphereCoveringInequality}
Let \(\Omega \subset \mathbb{R}^2\) be simply connected, and assume \(v_i \in C^2(\overline{\Omega})\), \(i = 1, 2\), satisfy
\begin{equation}\label{mainpde}
\Delta v_i + e^{2v_i} = f_i(y),
\end{equation}
where \(f_2 \ge f_1 \ge 0\) in \(\Omega\). Suppose that \(v_2 \ge v_1+c\), \(v_2 \not\equiv v_1\) in a piecewise Lipschitz subdomain \(\omega \subset \Omega\), and that \(v_2 - v_1 = c \geq 0\) on \(\partial \omega\). Then
\begin{equation}\label{GeneralSCI}
\int_{\omega} (e^{2v_1} + e^{2v_2})\,dy \;\ge\; 4(1 - \frac{c}{M})\pi,
\end{equation}
where
\[
 M = \max_{\omega}(v_2 - v_1).
\]
Moreover, equality holds if and only if \(c = 0\), \(f_2 \equiv f_1 \equiv 0\) in \(\omega\), and the weighted surfaces \((\omega, e^{2v_i} dy)\), \(i = 1, 2\), are isometric to two complementary spherical caps on the standard unit sphere.
\end{theorem}

Conceptually, this result may be viewed as a stability-type extension of the Sphere Covering Inequality, quantifying how much geometric information is lost when the boundary agreement condition is relaxed. Quantitative stability results for sharp conformally invariant inequalities have recently been obtained in different settings: Carlen~\cite{Carlen-2025-StabilityLogHLS} established stability for the logarithmic Hardy--Littlewood--Sobolev inequality as well as for the Onofri inequality, and Ghosh and Karmakar~\cite{GhoshKarmakar-arXiv2508.19930} demonstrated quantitative stability for the Chang--Gui inequality on the sphere.

\begin{remark}
The assumption \( v_2 \ge v_1 + c \) in Theorem \ref{GeneralSphereCoveringInequality} can be replaced by the weaker condition 
\[
M=\max_{\omega}(v_2 - v_1) > c.
\]
Indeed, under this weaker hypothesis one may apply Theorem~\ref{GeneralSphereCoveringInequality} to  the set 
\[
\omega' :=\{ x \in \omega : v_2(x) - v_1(x) > c \},
\]
and conclude that 
\[
\int_{\omega} (e^{2v_1} + e^{2v_2})\,dy 
\;\ge\; 
\int_{\omega'} (e^{2v_1} + e^{2v_2})\,dy 
\;\ge\; 
4\!\left(1 - \frac{c}{M}\right)\!\pi.
\]
\end{remark}

\section{Auxiliary Results and Preliminaries}
In this section, we establish several preliminary results that will be used in the proof of the main theorem. These lemmas, though elementary in nature, play an essential role in clarifying the analytic framework and preparing the ground for the more delicate geometric arguments that follow.

\begin{lemma}\label{CalcLemma0}
For \( k \in [0,1] \), define
\[
f(x) = \frac{x^{2-k} - x^{k}}{x^{2} - 1} + k - 1.
\]
Then \( f(x) \le 0 \) for all \( x \in (1,\infty) \).  
In particular, \( f(x) < 0 \) for all \( x \in (1,\infty) \) if \( k \in (0,1) \).
\end{lemma}

{\bf Proof.}
First observe that \( f \equiv 0 \) for \( k \in \{0,1\} \).  
Assume \( k \in (0,1) \). Since
\[
\lim_{x \to 1} f(x) = 0,
\]
it suffices to show that \( f \) is decreasing on \( (1,\infty) \).  
A direct computation gives
\[
f'(x) = \frac{x\big( -k x^{2-k} + (2-k)x^{k} - (2-k)x^{-k} + kx^{k-2} \big)}{(x^{2}-1)^{2}},
\]
so \( f \) is decreasing if
\[
f_1(x) := -k x^{2-k} + (2-k)x^{k} - (2-k)x^{-k} + kx^{k-2} < 0
\]
for all \( x > 1 \).  
Since \( f_1(1) = 0 \), it is enough to verify that \( f_1'(x) < 0 \) on \( (1,\infty) \), where
\[
f_1'(x) = k(2-k)\big(-x^{1-k} + x^{k-1} + x^{-k-1} - x^{k-3}\big).
\]
Note that \( f_1'(x) < 0 \) if
\[
f_2(x) := -x^{1-k} + x^{k-1} + x^{-k-1} - x^{k-3}
\]
is decreasing on \( (1,\infty) \).  
Indeed,
\begin{align*}
f_2'(x)
&= -(1-k)x^{-k} - (1-k)x^{k-2} - (k+1)x^{-k-2} + (3-k)x^{k-4} \\
&< \big( -(1-k) - (1-k) - (k+1) + (3-k) \big)x^{k-4} = 0,
\end{align*}
and hence \( f_2 \) is strictly decreasing on \( (1,\infty) \).  
The proof is complete.\hfill \( \Box \) \\ \\

We shall also need the following elementary result.

\begin{lemma}\label{CalcLemma}
Let \( a, b, k \) be positive constants with \( b > a > 0 \) and \( k \in [0,1] \). 
Suppose that
\begin{equation}\label{CalcLemmaAssumption}
\frac{8 + b^2x^2}{8 + a^2x^2} = \left( \frac{b}{a} \right)^k
\end{equation}
holds. Then the inequality
\begin{equation}\label{CalclemmaInequality}
\frac{a^2x^2}{8 + a^2x^2} + \frac{b^2x^2}{8 + b^2x^2} \ge k
\end{equation}
is satisfied. Moreover, if \( k \in (0,1) \), then the inequality is strict.
\end{lemma}

{\bf Proof.}
From \eqref{CalcLemmaAssumption}, we obtain
\begin{equation}
x^2 = \frac{8\left( \left( \frac{b}{a} \right)^k - 1 \right)}{b^2 - \left( \frac{b}{a} \right)^k a^2},
\end{equation}
and hence
\begin{equation}
8 + a^2x^2 = \frac{8 \left( \left( \frac{b}{a} \right)^2 - 1 \right)}{\left( \frac{b}{a} \right)^2 - \left( \frac{b}{a} \right)^k}.
\end{equation}
Consequently,
\begin{equation}
\frac{8}{8 + a^2x^2} = \frac{\left( \frac{b}{a} \right)^2 - \left( \frac{b}{a} \right)^k}{\left( \frac{b}{a} \right)^2 - 1}.
\end{equation}
Therefore,
\begin{align*}
\frac{a^2x^2}{8 + a^2x^2} + \frac{b^2x^2}{8 + b^2x^2}
&= 2 - 8 \left( \frac{1}{8 + a^2x^2} + \frac{1}{8 + b^2x^2} \right) \\
&= 2 - \frac{8}{8 + a^2x^2} \left( 1 + \frac{1}{\frac{8 + b^2x^2}{8 + a^2x^2}} \right) \\
&= 2 - \frac{\left( \frac{b}{a} \right)^2 - \left( \frac{b}{a} \right)^k}{\left( \frac{b}{a} \right)^2 - 1} 
      \left( 1 + \frac{1}{\left( \frac{b}{a} \right)^k} \right) \\
&= 2 - \frac{\left( \frac{b}{a} \right)^2 - \left( \frac{b}{a} \right)^k}{\left( \frac{b}{a} \right)^2 - 1} 
      \frac{1 + \left( \frac{b}{a} \right)^k}{\left( \frac{b}{a} \right)^k} \\
&= 2 - \frac{\left( \frac{b}{a} \right)^{2 - k} - 1}{\left( \frac{b}{a} \right)^2 - 1}
      \left( 1 + \left( \frac{b}{a} \right)^k \right) \\
&= 2 - \frac{\left( \frac{b}{a} \right)^{2 - k} - 1 + \left( \frac{b}{a} \right)^2 - \left( \frac{b}{a} \right)^k}
             {\left( \frac{b}{a} \right)^2 - 1} \\
&= 1 - \frac{\left( \frac{b}{a} \right)^{2 - k} - \left( \frac{b}{a} \right)^k}
             {\left( \frac{b}{a} \right)^2 - 1}.
\end{align*}
Hence, inequality \eqref{CalclemmaInequality} holds if and only if
\begin{equation}
\frac{\left( \frac{b}{a} \right)^{2 - k} - \left( \frac{b}{a} \right)^k}
     {\left( \frac{b}{a} \right)^2 - 1} \le 1 - k.
\end{equation}
By Lemma~\ref{CalcLemma0}, the above inequality holds for \( \frac{b}{a} = x > 1 \).\hfill \( \Box \) \\ \\

For \( \lambda>0 \), define \( U_\lambda \) by 
\begin{equation}\label{ULambda}
U_{\lambda}:=-2\ln(1+\frac{\lambda^2 |y|^2 }{8})+2\ln(\lambda),
\end{equation}
which satisfies the equation 
\[
\Delta U_\lambda +e^{U_\lambda}=0.\\
\]
\\
We are now ready to prove the main result of this section.  \\

\begin{theorem}\label{RadialTheorem}
Let \( b > a > 0 \), and let \( U_a \) and \( U_b \) be defined as in~\eqref{ULambda}. 
Let \( k \in (0,1] \), and suppose that 
\begin{equation}\label{UAssumption}
U_b - U_a = 2(1 - k)\ln\!\left( \frac{b}{a} \right)
\quad \text{on } \partial B_R,
\end{equation}
where \( B_R \) is the ball of radius \( R \) centered at the origin. Then 
\begin{equation}
\int_{B_R} \!\big(e^{U_a} + e^{U_b}\big)\, dy \ge 8k\pi,
\end{equation}
and the inequality is strict if \( k \in (0,1) \).
\end{theorem}

{\bf Proof.}
From \eqref{UAssumption}, we have 
\[
-2 \ln\!\left( \frac{8 + b^2R^2}{8 + a^2R^2} \right) 
  + 2\ln\!\left( \frac{b}{a} \right) 
  = 2(1 - k)\ln\!\left( \frac{b}{a} \right),
\]
and hence 
\[
\frac{8 + b^2R^2}{8 + a^2R^2} = \left( \frac{b}{a} \right)^k.
\]
It then follows from Lemma~\ref{CalcLemma} that 
\begin{align*}
\int_{B_R} \!\big(e^{U_a} + e^{U_b}\big)\, dy
&= 8\pi \left( 
      \frac{a^2R^2}{8 + a^2R^2} 
    + \frac{b^2R^2}{8 + b^2R^2} 
   \right)
   \ge 8k\pi,
\end{align*}
with strict inequality when \( k \in (0,1) \).\hfill \( \Box \) \\ \\ 

\begin{remark}
It is important to note that Lemmas~\ref{CalcLemma0}, ~\ref{CalcLemma}, and Theorem \ref{RadialTheorem} are not valid for \( k > 1 \). Consequently, Theorem~\ref{GeneralSphereCoveringInequality} does not hold in the case \( c < 0 \).\\
\end{remark}

\section{Generalized Sphere Covering Inequality}

In this section, we present the proof of our main result, Theorem~\ref{GeneralSphereCoveringInequality}.  As a first step, we derive a comparison principle for radial solutions  involved, which will play a crucial role in our proof.

\begin{proposition}\label{LastEstimate}
Assume that $\psi \in C^{0, 1} (\overline{B_R})$ is a strictly decreasing, radial, Lipschitz function,  and  satisfies 
\begin{equation}\label{superSol}
\int_{\partial B_r} |\nabla \psi|ds  \le \int_{B_r}e^{\psi}dy
\end{equation}
{\it a.e.}  $r\in (0,R)$ and $\psi(0)=U_{b}(0)$  for some $b>0$. Then 
\begin{equation}\label{lastEstimate}
  \int_{B_R} e^{\psi} dy\geq   \int_{B_R} e^{U_{b}} dy.
\end{equation}
Moreover if the inequality in (\ref{superSol}) is strict at some point in $(0,R)$, then the inequality in (\ref{lastEstimate}) is also strict.
\end{proposition}
{\bf Proof.} Let $0<b'<b$.  We claim that 
\begin{equation}\label{bPrimeIneq}
\int_{B_R}e^{\psi}dy> \int_{B_R}e^{U_{b'}}dy. 
\end{equation}

Since $\psi(0)=U_{b}(0)>U_{b'}(0)$ and $\psi$ is continuous, we must have $\psi>U_{b'}$ in a neighborhood of the origin.  If  $\psi \geq U_{b'}$ in $B_R$, then \eqref{lastEstimate} follows immediately.  Define 
\[r_1=\sup\{s\in [0,R): \psi (r) > U_{b'}(r) \ \ \text{for all} \ \ r\in [0,s) \}.\]
Then $0<r_1< R$, $\psi(r_1)=U_{b'}(r_1)$, and 
\begin{equation}\label{larger}
\int_{B_{r_1}}e^{\psi(y)}dy>\int_{B_{r_1}}e^{U_{b'}(y)}dy. 
\end{equation}
On the other hand, by Bol's inequality we have 
\begin{eqnarray*}
\frac{1}{2}\left( \int_{B_{r_1}} e^\psi\right) \left( 8 \pi-\int_{B_{r_1}}e^\psi \right)&\leq &\left( \int_{\partial B_{r_1}}e^{\frac{\psi}{2}}\right)^2= \left( \int_{\partial B_{r_1}}e^{\frac{U_{b'}}{2}}\right)^2\\
&=&\frac{1}{2}\left( \int_{B_{r_1}} e^{U_{b'}}\right) \left( 8 \pi-\int_{B_{r_1 }}e^{U_{b'}} \right),
\end{eqnarray*}
and therefore $\int_{B_{r_1}}e^{\psi}>4\pi$. This follows from\eqref{larger} and the fact that the polynomial $\frac{1}{2}x(8\pi-x)$ is strictly  increasing on $(0, 4\pi)$.  Now define 
\[r_2=\sup \{s\in [0,R]: \int_{B_{s}}e^{\psi(y)}dy\geq\int_{B_{s}}e^{U_{b'}(y)}dy\},\]
and note that $r_2 \geq r_1.$ We claim that $r_2=R$. Suppose this is not the case.  Then there exists $r_3$ with $R \geq r_3>r_2$ such that 
\[\psi(r_3) < U_{b'}(r_3) \ \ \text{and}\ \  \int_{B_{r_3}}e^{\psi(y)}dy< \int_{B_{r_3}}e^{U_{b'}(y)}dy \]
By Bol's inequality we have 
\begin{eqnarray*}
\frac{1}{2}\left( \int_{B_{r_3}} e^\psi\right) \left( 8 \pi-\int_{B_{r_3}}e^\psi \right)&\leq &\left( \int_{\partial B_{r_3}}e^{\frac{\psi}{2}}\right)^2< \left( \int_{\partial B_{r_3}}e^{\frac{U_{b'}}{2}}\right)^2\\
&=&\frac{1}{2}\left( \int_{B_{r_3}} e^{U_{b'}}\right) \left( 8 \pi-\int_{B_{r_3 }}e^{U_{b'}} \right),
\end{eqnarray*}
this is a contradiction because $\frac{1}{2}x(8\pi-x)$ is increasing only for $x\leq 4\pi $ and $\int_{B_{r_3}}e^{\psi}dy>4\pi$.  Hence the inequality  \eqref{bPrimeIneq} holds for all $0<b'<b$.  The inequality \eqref{lastEstimate} follows from \eqref{bPrimeIneq} by letting $b' \rightarrow b$. The proof is now complete. \hfill $\Box$ \\ \\

Having established the necessary analytic framework and comparison principles, we are now ready to prove the central theorem of this work, the generalized Sphere Covering Inequality.  
It is worth noting that Theorem~\ref{GeneralSphereCoveringInequality} follows directly from the next result after a simple change of variable \( w = 2v \). \\

\begin{theorem}
Let \(\Omega \subset \mathbb{R}^2\) be simply connected, and assume \(v_i \in C^2(\overline{\Omega})\), \(i = 1, 2\), satisfy
\begin{equation}\label{mainpde}
\Delta v_i + e^{v_i} = f_i(y),
\end{equation}
where \(f_2 \ge f_1 \ge 0\) in \(\Omega\). Suppose that \(v_2 \ge v_1+c\), \(v_2 \not\equiv v_1\) in a piecewise Lipschitz subdomain \(\omega \subset \Omega\), and that \(v_2 - v_1 = c \geq 0\) on \(\partial \omega\). Then
\begin{equation}\label{GeneralSCI}
\int_{\omega} (e^{v_1} + e^{v_2})\,dy \;\ge\; 8(1 - \frac{c}{M})\pi,
\end{equation}
where
\[
 M = \max_{\omega}(v_2 - v_1).
\]
Moreover, equality holds if and only if \(c = 0\), \(f_2 \equiv f_1 \equiv 0\) in \(\omega\), and the weighted surfaces \((\omega, e^{2v_i} dy)\), \(i = 1, 2\), are isometric to two complementary spherical caps on the standard unit sphere.
\end{theorem}

{\bf Proof.} Suppose $v_1$ and $v_2$ satisfy the assumptions of the theorem. Then 
\[\Delta (v_2-v_1)+e^{v_2}-e^{v_1}=f_2-f_1\geq 0.\]
We may assume $\int_{\omega}e^{v_1}dy<4\pi$ and 
\begin{equation}\label{a-choice}
    \int_{\omega} e^{v_1}+e^{v_2}dy<8\pi-\epsilon
\end{equation}
for some $\epsilon >0$.  For any $a>0$,  there exists $R_a>0$ such that 
\begin{equation}
\int_{\omega}e^{v_1}dy=\int_{B_{R_a}}e^{U_{a}}dy= \frac{8\pi a^2R_{a}^2}{8+a^2R_{a}^2},
\end{equation}
where $U_{a}$ be define by \eqref{ULambda}. Observe that $a^2R_a^2=K=\text{constant}$, for any choice of $a>0$. Now let $b=a+M$ where $M:=\max\limits_{\omega} (v_2 - v_1)$. Hence, we can choose $a$ large enough such that 
\begin{eqnarray*}
    0<\int_{B_{R_a}}e^{U_b}dy-\int_{B_{R_a}}e^{U_a}dy&=& 8\pi \left[ \frac{(a+M)^2R_{a}^2}{8+(a+M)^2R_{a}^2}- \frac{a^2R_{a}^2}{8+a^2R_{a}^2}\right]\\
    &=& 8\pi \left[ \frac{(a+M)^2\frac{K}{a^2}}{8+(a+M)^2\frac{K}{a^2}}- \frac{K}{8+K}\right]<\epsilon.
\end{eqnarray*}

Fix $a$ large enough such that 
\begin{equation}\label{aa-choice}
    \int_{B_{R_a}}e^{U_b}dy-\int_{B_{R_a}}e^{U_a}dy <\epsilon,
\end{equation}
and for $t>c$ let $ \omega_t =  \{y \in \omega:    v_2(y) -v_1 (y)> t\}$. Let $\phi$ be the symmetrization of $v_2-v_1$ with respect to the measures $e^{v_1}dy$ and $e^{U_{a}}dy$.  More precisely,  define
\[\omega_t:=\{x\in \omega: v_2(x)-v_1(x)>t\},\]
and define $\omega^*_t$ to be the ball centered at the  origin in $\R^2$ such that 
\[\int_{\omega^*_t}e^{U_{a}}dy=\int_{\omega_t}e^{v_1}dy:=\beta(t).\]

Then $\beta(t) $ is  a right-continuous function,  and $\phi: \omega \rightarrow \R$ defined by $\phi(y):=\sup \{t\in \R: y\in \omega^*_t\}$ provides an equimeasurable rearrangement of $v_2-v_1$ with respect to the measure $e^{v_1}dy$ and $e^{U_{a}}dy$, i.e. 
\begin{equation}\label{rearrang}
\int_{\{\phi>t\}}e^{U_{a}}dy=\int_{\{v_2-v_1>t\}}e^{v_1}dy, \ \ \forall  t>\min_{y\in \overline \omega}\phi,
\end{equation}
where $\omega^*$  a ball of radius $R_a<\infty$ and centered at the origin with $\int_{\omega^*}e^{U_{a}}dy=\int_{\omega}e^{v_1}dy$.

By Propositions 2.2 and 2.3 in \cite{GM},   Green's formula,   and Cavalieri's principle we have  
\begin{eqnarray} \label{firstIneq}
\int_{\{\phi=t\}}|\nabla \phi|ds &\leq &  \int_{\{v_2-v_1=t\}\cap \omega}|\nabla (v_2-v_1)|  ds    \\
&\leq&\int_{\omega_t} \bigl(e^{v_2}-e^{v_1} \bigr) dy   \quad \hbox{ (by Green's formula and equation \eqref{mainpde} ) }  \nonumber \\
&=& \int_{\omega_t}e^{v_2-v_1} e^{v_1}dy -\int_{\omega_t}e^{v_1}  dy \nonumber \\
&=& \int_{\{\phi>t\}}e^{\phi}e^{U_{a }}dy- \int_{\{\phi>t\}}e^{U_{a}} dy \nonumber \\
&&  \quad \hbox{ (by rearrangement and Cavalieri's principle) }  \nonumber \\
&=& \int_{\{\phi>t\}}e^{U_{a}+\phi}dy-\int_{\{\phi=t\}} |\nabla U_{a}|ds,  \quad  \hbox{for  {\it a.e.} }  t>c.\nonumber
\end{eqnarray} 
  Hence
\begin{equation}\label{Sub}
\int_{\{\phi=t\}} |\nabla (U_{a}+\phi)| ds\le \int_{\{\phi>t\} } e^{(U_{a}+\phi)}dy,  \quad  \hbox{for  {\it a.e.} }  t>c.
\end{equation}
  Since  $ \phi \geq c$ is decreasing in $r$,     $\psi:= U_{a}+\phi$ is a strictly decreasing radial function, and 
\begin{equation}\label{supersolution}
\int_{\partial B_r} |\nabla \psi|ds \le  \int_{B_r} e^{\psi} dy, \quad {\it a.e.}  \quad r \in (0, R_a), 
\end{equation}
by Proposition 2.3 in \cite{GM} and the above inequality we know that  $\psi$ belongs to  $C^{0, 1} (B_{R_a})$.  It follows from Proposition \ref{LastEstimate} that 
\[\int_{B_{R_a}}e^{\psi}dy =\int_{B_{R_a}}e^{U_{a}+\phi}dy \geq 
\int_{B_{R_a}}e^{U_{b}}dy.\]

Next we claim that $\psi(R_a) \geq U_{b} (R_a)$. Suppose that this is not the case. Since, by Proposition \ref{LastEstimate},
\[\int_{B_r}e^{\psi}dy \geq \int_{B_r}e^{U_b}dy, \ \ \forall r\in (0,R_a),\]
there must exist $R_0 \in (0,R_a)$ such that $\psi (R_0)=U_b(R_0)$, and  $\int_{B_{R_0}}e^{\psi}dy > \int_{B_{R_0}}e^{U_b}dy$. Thus it follows from the Sphere Covering Inequality that  
\[\int_{B_{R_0}}e^{\psi}+e^{U_b}dx \geq 8\pi.\]
Thus it follows from \eqref{aa-choice} that
\begin{eqnarray*}
    \int_{\omega}e^{v_1}+e^{v_2}dy &\geq& \int_{B_{R_0}}e^{U_a}+e^{\psi}dx\\
    &\geq & \int_{B_{R_0}}e^{U_b}+e^{\psi}dx-\epsilon \\
    &\geq & 8\pi -\epsilon. 
\end{eqnarray*}
This is a contradiction in view of \eqref{a-choice}, and hence $$\psi(R_a) \geq U_{b} (R_a)$$.  

Now we have
\begin{eqnarray*}
\int_{\Omega}(e^{v_1}+e^{v_2})dy&=& \int_{B_{R_a}}(e^{U_{a}}+e^{\psi})dy \\ &\geq& \int_{B_{R_a}}(e^{U_{a}}+e^{U_{b}})dy\geq 8\pi \left[  1-\frac{U_b(R_a)-U_a({R_a})}{M}\right] \\ 
&\geq& 8\pi \left[  1-\frac{\psi(R_a)-U_a({R_a})}{M}\right]\\
&= &8\pi \left[  1-\frac{(v_2-v_1)|_{\partial \omega}}{M}\right] = 8(1 - \frac{c}{M})\pi.
\end{eqnarray*}
Moreover,  if  the equality holds,  then  $k=1$, and the inequality \eqref{GeneralSCI} simplifies to the Sphere Covering Inequality \cite{GM} for which the equality case is established. \hfill $\Box$ \\

$\mathbf{Acknowledgement}.$ The first author is  partially supported by University of Macau research grants (No. CPG2024-00016-FST, CPG2025-00032-FST, SRG2023-00011-FST and MYRGGRG2023-00139-FST-UMDF), UMDF Professorial Fellowship of Mathematics, Macao SAR FDCT 0003/2023/RIA1 and Macao SAR FDCT 0024/2023/RIB1. The second author is supported in part by NSF grant DMS-1715850.

\bibliographystyle{plain}
\bibliography{GeneralizedSphereCoveringInequality.bib}

\end{document}